\documentclass[11pt,twoside]{amsart}

\usepackage{amsfonts}

\title{An axiomatic approach for degenerations in triangulated categories}

\author{Manuel Saorin}
\address{\newline
Departemento de Matem\'aticas,
\newline Universidad de Murcia, Aptdo. 4021
\newline 30100 Espinardo, Murcia,
\newline Spain}
\email{msaorinc@um.es}

\author{Alexander Zimmermann}
\address{\newline
Universit\'e de Picardie,
\newline D\'epartement de Math\'ematiques et LAMFA (UMR 7352 du CNRS),
\newline 33 rue St Leu,
\newline F-80039 Amiens Cedex 1,
\newline France}
\email{alexander.zimmermann@u-picardie.fr}

\newtheorem*{Theo2}{{Theorem}}
\newtheorem{Lemma1}{{Lemma}}
\newtheorem{Theo1}{{Theorem}}
\newtheorem{Def1}[Lemma1]{{Definition}}
\newtheorem{Prop1}[Lemma1]{{Proposition}}
\newtheorem{Claim1}[Lemma1]{{Claim}}
\newtheorem{Rem1}[Lemma1]{{Remark}}
\newtheorem{Cor1}[Lemma1]{{Corollary}}
\newtheorem{Ex1}[Lemma1]{{Example}}

\newenvironment{Lemma}{\begin{Lemma1}}{\end{Lemma1}}
\newenvironment{Def}{\begin{Def1}\em}{\end{Def1}}
\newenvironment{Prop}{\begin{Prop1}}{\end{Prop1}}
\newenvironment{Rem}{\begin{Rem1}\rm}{\end{Rem1}}
\newenvironment{Theorem}{\begin{Theo1}}{\end{Theo1}}

\newenvironment{Cor}{\begin{Cor1}}{\end{Cor1}}

\newenvironment{Example}{\begin{Ex1}\em}{\end{Ex1}}

\newcommand{\uar}{\uparrow}
\newcommand{\dar}{\downarrow}
\newcommand{\lra}{\longrightarrow}

\newcommand{\ra}{\rightarrow}
\newcommand{\sdp}{\times\kern-.2em\vrule height1.1ex depth-.05ex}
\newcommand{\epi}{\lra \kern-.8em\ra}

\newcommand{\N}{{\mathbb N}}

\newcommand{\Z}{{\mathbb Z}}

\newcommand{\cone}{\text{cone}}

 \normalbaselines
\subjclass[2010]{Primary 18E30; Secondary   16G10; 16G30; 16E35; 18E35 }

\date{May 20, 2014}
\thanks{Saor\'\i n is supported by research projects from the
Secretar\'\i a de Estado de Investigaci\'on,
Desarrollo e Innovaci\'on of the Spanish Government and the Fundaci\'on
'S\'eneca' of Murcia, with a part of FEDER funds.
Zimmermann is supported by STIC Asie project 'Escap' of the Minist\`ere
des Affaires \'Etrang\`eres de la France.}

\begin{document}

\begin{abstract}
We generalise Yoshino's definition of a degeneration
of two Cohen Macaulay modules to
a definition of degeneration between two objects in a
triangulated category. We derive some natural properties
for the triangulated category and the degeneration under
which the Yoshino-style degeneration is equivalent to the
degeneration defined by a specific distinguished triangle
analogous to Zwara's characterisation of degeneration in module varieties.
\end{abstract}

\maketitle

\section*{Introduction}

For an algebra $A$ over a field $k$, $A$-modules of fixed dimension $d$
over $k$ are $k$-vector spaces together with an action of $A$.
Fix $k^d$ and then
the conditions reflecting these properties can be formulated
in terms of vanishing of a finite number of polynomial equations.
Hence, the points of the algebraic variety defined by these polynomial
equations correspond to the possible $A$-module structures
on $k^d$. Two such structures are isomorphic if and only if they belong to
the same $Gl_d(k)$-orbit. The variety together with the
action of the linear group is called the module variety $mod_d(A)$.
If $M$ corresponds to the point $m$ and if $N$ corresponds
to the point $n$, then we say that $M$ degenerates to $N$ if
$n$ is in the Zariski-closure of the $Gl_d(k)$-orbit of $m$, and we denote
in this case $M\leq_{deg}N$.

Riedtmann and Zwara showed in \cite{Riedtmann,Zwara} that
$M\leq_{deg}N$ if and only if there is an $A$-module $Z$
and an exact sequence
$$0\lra Z\lra Z\oplus M\lra N\lra 0$$
and this is equivalent to the existence of an $A$-module $Z'$
and an exact sequence
$$0\lra N\lra Z'\oplus M\lra Z'\lra 0.$$
The concept of degeneration of modules is highly successful
in representation theory of algebras.

Various attempts were undertaken to prove an analogous result for
triangulated categories. In \cite{HuisgenSaorin} the first author
defined degenerations of complexes of $A$-modules, and a slight modification yields to
 a degeneration concept of the derived category. In \cite{JSZdegen} Jensen, Su and
the second author defined a topological space,
whose points correspond to objects in $D^-(A)$,
and on which a topological group acts. In a parallel development
Jensen, Madsen and Su give in \cite{JMS}
another algebraic variety with a group
action parameterising objects in the derived category of
$A_\infty$-modules over $A_\infty$-algebras.
Degeneration can be defined analogously to the module
variety situation and then it can be shown in the last two settings that
an object $M$ degenerates to an object $N$ if and only if
there is an object $Z$ and a distinguished triangle
$$N\lra Z\oplus M\lra Z\lra N[1].$$
This second condition on the existence of distinguished triangles
can be used to define some relation $\leq_\Delta$
on isomorphism classes of objects in triangulated categories $\mathcal T$.
It is shown in \cite{JSZpartord} that under certain conditions
the relation $\leq_\Delta$ is a partial order on the isomorphism classes of objects of $\mathcal T$.

As for stable categories Yoshino gave a definition of stable degeneration
$\leq_{stdeg}$ in \cite{Yoshinostable}
building on earlier work on an alternative approach for module varieties
(cf \cite{YoshinoCM,YoshinoModules}). For
stable categories of maximal Cohen-Macaulay modules over isolated singularities
Yoshino proved that
$\leq_{stdeg}$ is equivalent to the partial order
$\leq_\Delta$. Yoshino's construction uses that an object $M$
degenerates to $N$ if and only if there is a line through $N$
so that the
generic point of the line is isomorphic to the generic point
generated by $M$. Modelling this gives a criterion on the existence of
a special object over an extension of the base category.

A geometrically inspired notion of degeneration was still missing. The purpose of
this note is to develop such a notion which works in a very general
triangulated setting, and our main result Theorem~\ref{main} shows
that this geometrically inspired notion of degeneration is
equivalent to $\leq_\Delta$. As a by-product the notion also works
for bounded derived categories improving in this way
\cite{HuisgenSaorin}, \cite{JSZdegen} and \cite{JMS}. We proceed
analogously to Yoshino's approach and formalise the existence of a
line to a functor to an extension category together with an element
in its centre, and the condition of the generic point to a condition
on the localisation category. The main definition is
Definition~\ref{degendatadef} and the definition of a degeneration
is given in Definition~\ref{degendef}.

We give an outline of the structure of the paper.
In Section~\ref{themaindefinitionsection} the definition of
triangulated degeneration is given.
In Section~\ref{triangulatedtoolssection} we
give some technical details valid for triangulated categories
which we use in the sequel and which we recall for the reader's convenience.
Our main result is then formulated and proved in
Section~\ref{mainsection}, which contains
Theorem~\ref{main} and its proof fills  the entire section.
We conclude with Section~\ref{conclusionsection}
which contains some consequences and remarks on degeneration
for triangulated categories.

\section{Defining degeneration}
\label{themaindefinitionsection}

\subsection{The ring theoretic degeneration axioms}

We shall define a system of axioms modelling the classical situation.
First recall Yoshino's definition of module degeneration.

\begin{Def} (Yoshino~\cite{YoshinoModules})
Let $k$ be a field and let $A$ be a $k$-algebra. Then for all finitely generated
$A$-modules $M$ and $N$ we say that $M$ degenerates to $N$ along a discrete
valuation ring if there is a discrete valuation ring $(V,tV,k)$ that is a
$k$-algebra (where $t$ is a prime element) and a finitely generated
$A\otimes_kV$-module $Q$ such that
\begin{enumerate}
\item
$Q$ is flat as a $V$-module
\item
$Q/tQ\simeq N$ as an $A$-module
\item $Q[\frac 1t]\simeq M\otimes_kV[\frac 1t]$ as an $A\otimes V[\frac 1t]$-module.
\end{enumerate}
\end{Def}

Yoshino shows in \cite{YoshinoModules} that $M$ degenerates to $N$ along a
discrete valuation ring if and only if there is a short exact sequence of
finitely generated $A$-modules
$$0\ra Z\stackrel{\begin{pmatrix}v\\ u\end{pmatrix}}{\lra}Z\oplus M\lra N\ra 0$$
such that $v$ is a nilpotent endomorphism of $Z$.

Yoshino generalizes in \cite{Yoshinostable} this concept in
the obvious way to define degeneration
in the (triangulated !) stable category of maximal Cohen-Macaulay modules over a
commutative local Gorenstein $k$-algebra.

\medskip

We transpose this to general triangulated categories in a
categorical framework. Observe first the following remark.

\begin{Rem}
Let $\mathcal{C}$ be a triangulated category and let $t$ be an
element of its center, i.e. a natural transformation of triangulated
functors $t:\textrm{id}_\mathcal{C}\longrightarrow \textrm{id}_\mathcal{C}$, where
$\textrm{id}_\mathcal{C}$ is the identity functor on $\mathcal C$. Then $\Sigma:=\{t_X^n|\;X
\mbox{ object of }{\mathcal C},\;n\in\N\}$ is, in the terminology of
Verdier (see \cite[Section II.2]{Verdier}), a multiplicative system
of $\mathcal{C}$ compatible with the triangulation. The localization
of $\mathcal{C}$ with respect to $\Sigma$, in the sense of
Gabriel-Zisman~\cite{GabrielZisman},  is denoted by $\mathcal{C}[t^{-1}]$ in the
sequel and is then a triangulated category such that the canonical
functor $p:\mathcal{C}\longrightarrow\mathcal{C}[t^{-1}]$ is
triangulated. Note that the Hom spaces in $\mathcal{C}[t^{-1}]$ are
sets, so that this is a proper category. If $\mathcal{D}$ is any
triangulated category, then giving a triangulated functor $\Psi
:\mathcal{C}[t^{-1}]\longrightarrow\mathcal{D}$ is equivalent to
giving a triangulated functor, which we denote the same,  $\Psi
:\mathcal{C}\longrightarrow\mathcal{D}$ such that $\Psi (t_X)$ is an
isomorphism.
\end{Rem}

\begin{Def}\label{degendatadef}
Let $k$ be a commutative ring and let ${\mathcal C}_k^\circ$ be a
$k$-linear triangulated category with split idempotents.

A degeneration data for ${\mathcal C}_k^\circ$ is given by
\begin{itemize}
\item a triangulated category ${\mathcal C}_k$ with split idempotents and a fully faithful
embedding ${\mathcal C}_k^\circ\longrightarrow{\mathcal C}_k$,
\item
a triangulated category   ${\mathcal C}_V$ with split idempotents
and a full triangulated subcategory ${\mathcal C}_V^\circ$,
\item  triangulated functors $\uar_k^V:{\mathcal C}_k\lra {\mathcal C}_V$
and  $\Phi:{\mathcal C}_V^\circ\ra {\mathcal C}_k$, so that
$({\mathcal C}_k^\circ)\uar_k^V\subseteq {\mathcal C}_V^\circ$, when
we view ${\mathcal C}_k^\circ$ as a full subcategory of ${\mathcal
C}_k$,
\item  a natural transformation
$\text{id}_{{\mathcal C}_V}\stackrel{t}{\lra} \text{id}_{{\mathcal C}_V}$
of triangulated functors
\end{itemize}
These triangulated categories and functors should satisfy the following axioms:
\begin{enumerate}
\item \label{4}
For each object $M$ of ${\mathcal C}_k^\circ$ the morphism
$\Phi(M\uar_k^V)\stackrel{\Phi(t_{M\uar_k^V})}{\lra}\Phi(M\uar_k^V)$
is a split monomorphism in ${\mathcal C}_k$.
\item \label{6}
For all objects $M$ of ${\mathcal C}_k^\circ$ we get
$\Phi(\cone(t_{M\uar_k^V}))\simeq M$.
\end{enumerate}
\end{Def}

All throughout the paper, whenever we have a degeneration data for
$\mathcal{C}_k^\circ$ as above, we will see $\mathcal{C}_k^\circ$ as
a full subcategory of $\mathcal{C}_k$.

\begin{Def}\label{degendef}
Given two objects $M$ and $N$ of ${\mathcal C}_k^\circ$ we say that
$M$ degenerates to $N$ in the categorical sense if there is a
degeneration data for ${\mathcal C}_k^\circ$ and an object $Q$ of
${\mathcal C}_V^\circ$ so that
$$p(Q)\simeq p(M\uar_k^V)\mbox{ in ${\mathcal C}_V^\circ[t^{-1}]$
and }\Phi(\cone(t_Q))\simeq N,$$ where $p:{\mathcal
C}_V^\circ\longrightarrow{\mathcal C}_V^\circ [t^{-1}]$ is the
canonical functor. In this case we write $M\leq_{cdeg}N$.
\end{Def}

\section{Some tools in triangulated categories}

\label{triangulatedtoolssection}

We will use a well-known result in triangulated categories, which can be found
in e.g. \cite{May} or \cite[Lemma 3.4.5]{reptheobuch}.

\begin{Lemma}\label{Mayslemma}
Let $\mathcal T$ be a triangulated category and let
$$\begin{array}{ccc}
A&\stackrel{\alpha}{\ra}&B\\
\phantom{\gamma}\dar\gamma&&\phantom{\beta}\dar\beta\\
C&\stackrel{\delta}{\ra}&D
\end{array}$$
be a commutative diagram in $\mathcal T$. Then
we may complete the horizontal and vertical maps to distinguished triangles
and there are morphisms $\epsilon$ and $\varphi$ so that in the diagram
$$\begin{array}{cccccccccc}
A&\stackrel{\alpha}{\ra}&B&\ra&C(\alpha)&\ra&A[1]\\
\phantom{\gamma}\dar\gamma&&\phantom{\beta}\dar\beta&&
\phantom{\epsilon}\dar\epsilon&&\phantom{\gamma[1]}\dar\gamma[1]\\
C&\stackrel{\delta}{\ra}&D&\ra&C(\delta)&\ra&C[1]\\
\dar&&\dar&&\dar&&\dar\\
C(\gamma)&\ra&C(\beta)&\ra&C(\epsilon)&\ra&C(\gamma)[1]\\
\dar&&\dar&&\dar&\mbox{\scriptsize ``$(-)$''}&\dar\\
A[1]&\stackrel{\alpha[1]}{\ra}&B[1]&\ra&C(\alpha)[1]&\ra&A[2]\\
\end{array}$$
all horizontal and vertical sequences are distinguished triangles.
All squares are commutative, except the lower right one, which is anticommutative.
\end{Lemma}

Idea of the proof.
The construction of $\epsilon$ is the following.
Take the composition $\xi:=\beta\circ\alpha=\delta\circ\gamma$
and consider the distinguished triangle
$$A\stackrel{\xi}{\ra}D\ra V\ra A[1].$$
Then the octahedral axiom applied to the two factorisations of $\xi$
give distinguished triangles
$$C(\alpha)\stackrel{\rho}{\lra}V\lra C(\beta)\ra C(\alpha)[1]$$
and
$$C(\gamma)\ra V\stackrel{\sigma}{\lra}C(\delta)\ra C(\gamma)[1].$$
Define $\epsilon:=\sigma\circ\rho$.

\begin{Rem}\label{octrahedralmapcanbechosen}
Recall the octahedral axiom.
Given three objects $X_1,X_2,X_3$ and suppose $\alpha_2:X_1\lra X_3$
factorises $\alpha_2=\alpha_1\circ\alpha_3$ for $\alpha_3\in \text{Mor}_{\mathcal T}(X_1,X_2)$ and
$\alpha_1\in \text{Mor}_{\mathcal T}(X_2,X_3)$. Then forming the triangles
$$X_2\stackrel{\alpha_1}{\lra}X_3\stackrel{\beta_1}{\lra}Z_1\stackrel{\gamma_1}{\lra}X_2[1]$$
$$X_1\stackrel{\alpha_3}{\lra}X_2\stackrel{\beta_3}{\lra}Z_3\stackrel{\gamma_3}{\lra}X_1[1]$$
$$X_1\stackrel{\alpha_2}{\lra}X_3\stackrel{\beta_2}{\lra}Z_2\stackrel{\gamma_2}{\lra}X_1[1]$$
above $\alpha_1$, $\alpha_2$ and $\alpha_3$, there are morphisms
$\delta_1:Z_3\ra Z_2$, $\delta_3:Z_2\ra Z_1$ and $\delta_2:Z_1\ra Z_3[1]$
so that
$$Z_3\stackrel{\delta_1}{\lra}Z_2\stackrel{\delta_3}{\lra}Z_1\stackrel{\delta_2}{\lra}Z_3[1]$$
is a distinguished triangle, and so that
$$\gamma_2\circ\delta_1=\gamma_3\;;\;\delta_3\circ\beta_2=\beta_1\;;\;
\delta_2=\beta_3[1]\circ\gamma_1\;;$$
$$\beta_2\circ\alpha_1=\delta_1\circ\beta_3\mbox{
and }
\gamma_1\circ\delta_3=\alpha_3[1]\circ\gamma_2$$
as illustrated in the following picture.

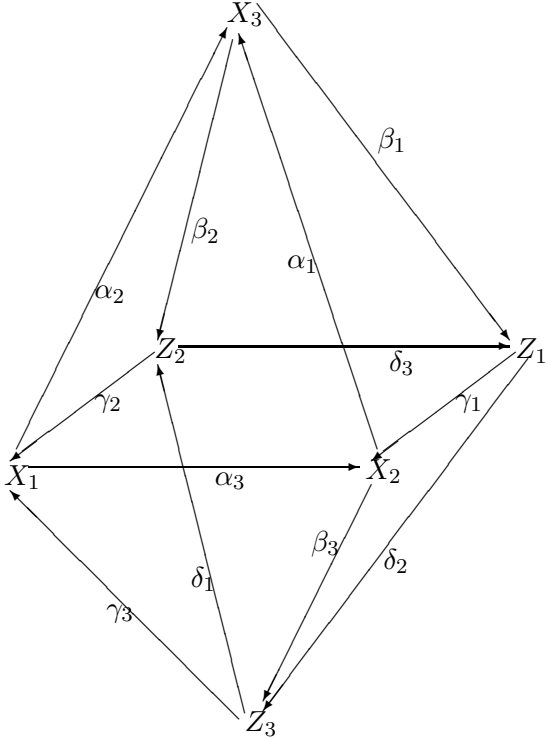
\unitlength.8cm
\begin{picture}(11,13.5)
\put(0,4.9){$X_1$}
\put(6,5){$X_2$}
\put(2.5,7){$Z_2$}
\put(8.5,7){$Z_1$}
\put(4,.8){$Z_3$}
\put(3.7,12.6){$X_3$}
\put(.4,5.2){\vector(1,0){5.5}}\put(3.5,4.9){$\alpha_3$}
\put(2.9,7.2){\vector(1,0){5.5}}\put(6.4,6.8){$\delta_3$}
\put(2.5,7.1){\vector(-4,-3){2.4}}\put(1.5,6.2){$\gamma_2$}
\put(8.5,7.1){\vector(-4,-3){2.4}}\put(7.5,6.2){$\gamma_1$}
\put(4,1.1){\vector(-1,4){1.45}}\put(3.1,3.2){$\delta_1$}
\put(3.9,1){\vector(-1,1){3.8}}\put(1.7,2.7){$\gamma_3$}
\put(6.1,4.9){\vector(-1,-2){1.8}}\put(5.1,3.8){$\beta_3$}
\put(8.7,7){\vector(-3,-4){4.4}}\put(6.3,3.5){$\delta_2$}
\put(.2,5.5){\vector(1,2){3.5}}\put(1.5,8){$\alpha_2$}
\put(3.8,12.3){\vector(-1,-4){1.25}}\put(3.1,9){$\beta_2$}
\put(6.2,5.5){\vector(-1,3){2.3}}\put(4.7,8.5){$\alpha_1$}
\put(4.2,12.9){\vector(3,-4){4.2}}\put(6.2,10.5){$\beta_1$}
\end{picture}

We remind the reader that the morphisms denoted by $\gamma$ are all of degree $1$,
and that also $\delta_2$ is of degree $1$.

Recall that the octahedral axiom just guarantees the existence of the
distinguished triangle as claimed. No uniqueness is claimed.
However, if $\gamma_3=0$ then
$\beta_3$ is a split epimorphism, so that there is a morphism
$Z_3\stackrel{\omega_3}{\lra}X_2$ with $\beta_3\circ\omega_3=\text{id}_{Z_3}$.
Since $\delta_1$
satisfies $\beta_2\circ\alpha_1=\delta_1\circ\beta_3$, and since split epimorphisms
are epimorphisms, this equation determines $\delta_1$.
More precisely, multiply with
$\omega_3$ from the right to obtain
$\beta_2\circ\alpha_1\circ\omega_3=\delta_1\circ\beta_3\circ\omega_3=\delta_1.$
\end{Rem}

We shall need a technical and well-known lemma on Gabriel-Zisman localisations.
For the convenience of the reader we include the short proof.

\begin{Lemma}\label{LemmaB}
Let $X$ be an object of ${\mathcal C}_V^\circ$.
We denote by $p:{\mathcal C}_V^\circ\lra {\mathcal C}_V^\circ[t^{-1}]$
the natural functor for Gabriel-Zisman localisation.
Then the following assertions are equivalent.
\begin{itemize}
\item $p(X)=0$
\item there is $n\geq 0$ so that $t_X^n=0$.
\end{itemize}
\end{Lemma}

\begin{proof}
Put again $\Sigma:=\{t_X^n|\;X \mbox{ object of }{\mathcal C}_1,\;n\in\N\}$.
A morphism $f$ of ${\mathcal C}_V^\circ$ has the property $p(f)=0$
if and only if there is an $s\in\Sigma$ so that $f\circ s=0$. Hence
\begin{eqnarray*}
p(X)=0&\Leftrightarrow& p(\text{id}_X)=0\\
&\Leftrightarrow& \exists s\in \Sigma:\text{id}_X\circ s=0\\
&\Leftrightarrow& \mbox{ the zero endomorphism on $X$ is in $\Sigma$.}
\end{eqnarray*}
Hence, $$p(X)=0\Leftrightarrow\exists n\in \N: t_X^n=0.$$
\end{proof}

\section{The theorem and its proof}
\label{mainsection}

\subsection{Categorical degeneration implies triangle degeneration}

\begin{Prop}\label{catimpliestriangle}
Let ${\mathcal C}_k^\circ$ be a triangulated $k$-category with split
idempotents and let $M$ and $N$ be two objects of ${\mathcal
C}_k^\circ$. Then $M\leq_{cdeg}N$ implies that there is  a
distinguished triangle
$$Z\stackrel{\begin{pmatrix}v\\ u\end{pmatrix}}{\lra} Z\oplus M\lra N\lra Z[1]$$
in $\mathcal{C}_k^\circ$ with nilpotent endomorphism $v$ of $Z$.
\end{Prop}

\begin{proof} Suppose  that $M\leq_{cdeg}N$, and there is hence a
degeneration data $({\mathcal C}_k,{\mathcal C}_V,\uar_k^V,t),$ an
object $Q$ of ${\mathcal C}_V^\circ$, so that $p(Q)\simeq
p(M\uar_k^V)$ in ${\mathcal C}_V^\circ[t^{-1}]$. Therefore there is a
morphism $M\uar_k^V\stackrel{f}{\lra}Q$ and a morphism
$Q\stackrel{g}{\lra}M\uar_k^V$ in ${\mathcal C}_V^\circ$ such that
$g\circ f\circ t_{M\uar_k^V}^r=t_{M\uar_k^V}^s$ and $f\circ g\circ
t_Q^m=t_Q^n$, for some  $r,s,m,n\in\mathbb{N}$. From  the fact that
$t$ is a natural transformation, we also get that
$t_{M\uar_k^V}^r\circ g\circ f=t_{M\uar_k^V}^s$ and, by applying the
functor $\Phi:{\mathcal C}_V^\circ\lra {\mathcal C}_k$, we get that
 $\Phi (f)$ is a split monomorphism.

By completing $f$ to a triangle in $\mathcal{C}_V$

$$M\uar_k^V\stackrel{f}{\lra}Q\stackrel{\varphi}{\lra}
\cone(f)\lra M\uar_k^V[1],$$  we then get a split distinguished
triangle in $\mathcal{C}_k$
$$\Phi(M\uar_k^V)\stackrel{\Phi(f)}{\lra}\Phi(Q)\stackrel{\Phi(\varphi)}{\lra}
\Phi(\cone(f))\stackrel{0}{\lra}\Phi(M\uar_k^V)[1].$$

By axiom~(\ref{4}) we get a
distinguished triangle
$$\Phi(M\uar_k^V)\stackrel{\Phi(t_{M\uar_k^V})}{\lra}
\Phi(M\uar_k^V)\stackrel{\lambda}{\lra}M\stackrel{0}{\lra}\Phi(M\uar_k^V)[1].$$
We define $\mu:=\Phi(f)\circ\Phi(t_{M\uar_k^V})$ and apply the octahedral axiom
to this factorisation.
Since $\Phi$ is assumed to be triangulated we get that $\Phi(\cone(f))=\cone(\Phi(f))$
and we therefore obtain the commutative octahedral axiom diagram
\unitlength1cm
\begin{center}
\begin{picture}(12,7)
\put(-1.5,3){$(I)$}
\put(-.3,6.3){$\Phi(M\uar_k^V)$}
\put(.8,4){$\Phi(M\uar_k^V)$}
\put(2,2){$M$}
\put(2.6,0){$\Phi(M\uar_k^V)[1]$}
\put(0,6.1){\vector(1,-2){.8}}
\put(1.2,3.8){\vector(1,-2){.7}}
\put(2.2,1.8){\vector(1,-2){.7}}
\put(-1,5){$\Phi(t_{M\uar_k^V})$}
\put(1.3,3){$\lambda$}
\put(2.3,1){$0$}
\put(2.5,4.2){\vector(1,0){1.4}}
\put(2.8,3.8){$\Phi(f)$}
\put(4,4){$\Phi(Q)$}
\put(5,4.2){\vector(1,0){1.4}}
\put(5.5,3.9){$\Phi(\varphi)$}
\put(6.6,4){$\Phi(\cone(f))$}
\put(8.5,4.2){\vector(1,0){1.4}}
\put(9,3.8){$0$}
\put(10.1,4){$\Phi(M\uar_k^V)[1]$}
\put(.5,6.1){\vector(2,-1){3.5}}
\put(2,5){$\mu$}
\put(4.9,3.8){\vector(2,-1){1.5}}
\put(5,3.4){$\tau$}
\put(6,2.8){$\cone(\mu)$}
\put(6.9,2.7){\vector(2,-1){1.5}}
\put(8.5,1.7){$\Phi(M\uar_k^V)[1]$}
\put(2.4,2.1){\vector(4,1){3.5}}
\put(6.5,3.2){\vector(1,2){.35}}
\put(6.8,3.5){$\psi$}
\put(6.9,4.4){\vector(0,1){1.3}}
\put(6.6,6){$M[1]$}
\end{picture}
\end{center}
where
$$M\lra \cone(\mu)\lra\Phi(\cone(f))\lra M[1]$$
is a distinguished triangle.
Since $\Phi(f)$ is a split monomorphism we
further get that
$\Phi(\varphi)$ is a split epimorphism with right inverse denoted
$\omega:\Phi(\cone(f))\lra \Phi(Q)$, i.e.
$\Phi(\varphi)\circ\omega=\textrm{id}_{\Phi(\cone(f))}$.
But then
$$\psi\circ\tau\circ\omega=\Phi(\varphi)\circ\omega=\text{id}_{\Phi(\cone(f))}$$
and therefore $\psi$ is a split epimorphism.
This shows that $$\cone(\mu)\simeq\Phi(\cone(f))\oplus M.$$

Hence
$$\Phi(\cone(f))\oplus M\simeq \cone(\mu)\simeq \cone(\Phi(f\circ t_{M\uar_k^V})).$$
Observe that, since $t$ is a natural transformation, we get
$$t_Q\circ f=f\circ t_{M\uar_k^V}.$$

The definition of categorical degeneration implies that we have a
distinguished triangle
$$\Phi(Q)\stackrel{\Phi(t_Q)}{\lra}\Phi(Q)\lra N\lra \Phi(Q)[1].$$
Now we consider the factorisation
$$\mu=\Phi(t_Q\circ f)=\Phi(t_Q)\circ\Phi(f)$$
and obtain by the octahedral axiom the following diagram.

\unitlength1cm
\begin{center}
\begin{picture}(12,6)
\put(-1.5,3){$(II)$}
\put(0,5){$\Phi(M\uar_k^V)$}
\put(1.6,5.2){\vector(1,0){2.5}}
\put(2.4,5.3){$\Phi(f)$}
\put(2.6,4){$\mu$}
\put(4.4,5){$\Phi(Q)$}
\put(5.3,5.2){\vector(1,0){2.5}}
\put(6,5.3){$\Phi(\varphi)$}
\put(8,5){$\Phi(\cone(f))$}
\put(10,5.2){\vector(1,0){1}}
\put(10.5,5.3){$0$}
\put(11.1,5){$\Phi(M\uar_k^V)[1]$}
\put(4.8,4.8){\vector(0,-1){1}}
\put(4.9,4.2){$\Phi(t_Q)$}
\put(4.4,3.4){$\Phi(Q)$}
\put(5.8,3.35){$\tau$}
\put(4.8,3.2){\vector(0,-1){1}}
\put(4.6,1.8){$N$}
\put(4.8,1.7){\vector(0,-1){1}}
\put(4.6,.4){$\Phi(Q)[1]$}
\put(1,4.8){\vector(3,-1){3.4}}
\put(5.3,3.4){\vector(3,-1){1.6}}
\put(6.1,2.6){$\Phi(\cone(t_Q\circ f))$}
\put(7.8,4){$\zeta$}
\put(8,2.5){\vector(3,-1){1.5}}
\put(8.5,2){$\nu$}
\put(9.5,1.6){$\Phi(M\uar_k^V)[1]$}
\put(8.5,4.9){\vector(-1,-2){1}}
\put(6.7,2.5){\vector(-4,-1){1.7}}
\put(4.5,1.99){\vector(-1,0){1.5}}
\put(.7,1.9){$\Phi(\cone(f))[1]$}
\end{picture}
\end{center}
where we thus obtain a distinguished triangle
$$\Phi(\cone(f))\stackrel{\zeta}{\lra}\Phi(\cone(t_Q\circ f))\lra N\lra \Phi(\cone(f))[1].$$
Since $\Phi(\varphi)$ is a split epimorphism, Remark~\ref{octrahedralmapcanbechosen}
shows that there is a unique choice
$\zeta:=\tau\circ\Phi(t_Q)\circ\omega$ so that the appearing diagrams
are all commutative.

However, we already identified
$$\Phi(\cone(t_Q\circ f))=\Phi(\cone(f\circ t_{M\uar_k^V}))\simeq
\Phi(\cone(f))\oplus M$$ so that we may put $Z:=\Phi(\cone(f))$ to
obtain a distinguished triangle in $\mathcal{C}_k$
$$Z\stackrel{\begin{pmatrix}v\\ u\end{pmatrix}}{\lra} Z\oplus M\stackrel{\begin{pmatrix}\ell&m\end{pmatrix}}{\lra} N\lra Z[1].$$

\medskip

We  need to show that $v$ is nilpotent.

Recall that $Z=\Phi(\cone(f))=\cone(\Phi(f))$.
The octahedral axiom
diagram $(I)$ implied the fact that
$$\Phi(\cone(t_Q\circ f))\simeq M\oplus\Phi(\cone(f))$$ and
the octahedral axiom diagram $(II)$ gives the distinguished triangle
$$\Phi(\cone(f))\stackrel{\zeta}{\lra}\Phi(\cone(t_Q\circ f))\lra N\lra \Phi(\cone(f))[1].$$
Since by the octahedral diagram $(I)$ we get
a split distinguished triangle
$$M\lra \Phi(\cone(t_Q\circ f))\stackrel{\psi}{\lra} Z\lra M[1]$$
the canonical projection $\Phi(\cone(t_Q\circ f))=Z\oplus M\lra Z$
is obtained by $\psi$, and hence
\begin{eqnarray*}
v&=&\psi\circ\zeta\\
&=&\psi\circ\tau\circ\Phi(t_Q)\circ\omega\\
&=&\Phi(\varphi)\circ\Phi(t_Q)\circ\omega\\
&=&\Phi(\varphi\circ t_Q)\circ\omega\\
&=&\Phi(t_{\cone(f)}\circ \varphi)\circ\omega\\
&=&\Phi(t_{\cone(f)})\circ \Phi(\varphi)\circ\omega\\
&=&\Phi(t_{\cone(f)})
\end{eqnarray*}
where we used that $\omega$ is right inverse to $\Phi(\varphi)$, and that
$t$ is a natural transformation. But since
$p(f)$ is an isomorphism in ${\mathcal C}_V^\circ[t^{-1}]$, we have that
$$p(\cone(f))\simeq \cone(p(f))\simeq 0$$
and we get that $t_{\cone(f)}$ is nilpotent by Lemma~\ref{LemmaB}. Therefore
$v=\Phi(t_{\cone(f)})$ is nilpotent.

\medskip

We need to show that $Z$ is in ${\mathcal C}_k^\circ$.
We shall apply here Lemma~\ref{Mayslemma} to the case
$$\begin{array}{ccc}
\Phi(M\uar_k^V)&\stackrel{\Phi(f)}{\ra}&\Phi(Q)\\
\phantom{\Phi(t_{M\uar_k^V})}\dar\Phi(t_{M\uar_k^V})&&\phantom{\Phi(t_Q)}\dar\Phi(t_Q)\\
\Phi(M\uar_k^V)&\stackrel{\Phi(f)}{\ra}&\Phi(Q)
\end{array}$$
Recall that in Lemma~\ref{Mayslemma} we apply the octahedral axiom to the
composition $$\Phi(t_Q)\circ\Phi(f)=\mu=\Phi(f)\circ\Phi(t_{M\uar_k^V}).$$
The decomposition
$\Phi(t_Q)\circ\Phi(f)=\mu$ gives the octahedral diagram $(II)$ and the decomposition
$\mu=\Phi(f)\circ\Phi(t_{M\uar_k^V})$ gives the octahedral diagram $(I)$.
The distinguished triangles we obtain from these diagrams
were determined in the first two steps of the proof.
For the diagram $(II)$ we get the distinguished triangle
$$Z\stackrel{\begin{pmatrix}v\\ u\end{pmatrix}}{\lra}Z\oplus M\lra N\lra Z[1]$$
and for the diagram $(I)$ we get the split distinguished triangle
$$M\lra Z\oplus M\stackrel{\begin{pmatrix}\text{id}_Z&0\end{pmatrix}}{\lra}Z\lra M[1].$$
We see that the map $\varepsilon$ in
$$\begin{array}{cccccccccc}
\Phi(M\uar_k^V)&\stackrel{\Phi(f)}{\ra}&\Phi(Q)&\ra&Z&\ra&\Phi(M\uar_k^V)[1]\\
\phantom{\Phi(t_{M\uar_k^V})}\dar\Phi(t_{M\uar_k^V})&&\phantom{\Phi(t_Q)}\dar\Phi(t_Q)&&
\phantom{\varepsilon}\dar\varepsilon&&\Phi(t_{M\uar_k^V})[1]\dar\phantom{\Phi(t_{M\uar_k^V})[1]}\\
\Phi(M\uar_k^V)&\stackrel{\Phi(f)}{\ra}&\Phi(Q)&\ra &Z&\ra&\Phi(M\uar_k^V)[1]\\
\dar&&\dar&&\dar&&\dar\\
M&\ra&N&\ra&C(\epsilon)&\ra&M[1]\\
\dar&&\dar&&\dar&\mbox{\scriptsize ``$(-)$''}&\dar\\
\Phi(M\uar_k^V)[1]&\stackrel{\Phi(f)[1]}{\ra}&\Phi(Q)[1]&\ra&Z[1]&\ra&\Phi(M\uar_k^V)[2],
\end{array}$$
such that all rows and columns are distinguished triangles, is obtained as
$$\varepsilon=\begin{pmatrix}\text{id}_Z&0\end{pmatrix}\cdot\begin{pmatrix}v\\ u\end{pmatrix}=v=\Phi(t_{\cone(f)}).$$
We have already seen  that $t_{\cone(f)}$ is nilpotent, and so is therefore $v$.
We remark that the left most column is a split triangle, and hence there is a unique
choice for the map $M\lra N$, as can be obtained similarly as in
Remark~\ref{octrahedralmapcanbechosen}.
We now use the second lower row distinguished triangle
$$M\ra N\ra \cone(\Phi(t_{\cone(f)}))\ra M[1].$$
Since $M$ and $N$ are in the triangulated category ${\mathcal C}_k^\circ$, and since
$$M\ra N\ra \cone(\Phi(t_{\cone(f)}))\ra M[1]$$
is a distinguished triangle, also $\cone(\Phi(t_{\cone(f)}))$ is in ${\mathcal C}_k^\circ$.
We apply the octahedral axiom to the factorisation
$$\Phi(t_{\cone(f)}^2)=\Phi(t_{\cone(f)})\circ \Phi(t_{\cone(f)})$$
and obtain a distinguished triangle
$$\cone(\Phi(t_{\cone(f)}))\ra \cone(\Phi(t_{\cone(f)})^2)\ra \cone(\Phi(t_{\cone(f)}))\ra
\cone(\Phi(t_{\cone(f)}))[1]$$
and hence also $\cone(\Phi(t_{\cone(f)})^2)$ is in ${\mathcal C}_k^\circ$.
The octahedral axiom applied to
$$\Phi(t_{\cone(f)}^n)=\Phi(t_{\cone(f)}^{n-1})\circ \Phi(t_{\cone(f)})$$
yields distinguished triangles
$$\cone(\Phi(t_{\cone(f)}))\ra \cone(\Phi(t_{\cone(f)})^n)\ra \cone(\Phi(t_{\cone(f)})^{n-1})\ra
\cone(\Phi(t_{\cone(f)}))[1]$$
and by induction we get $\cone(\Phi(t_{\cone(f)})^n)$ is in ${\mathcal C}_k^\circ$
for all $n\geq 1$.
But $\Phi(t_{\cone(f)})$ is nilpotent, and so
the distinguished triangle
$$\Phi(\cone(f))\stackrel{\Phi(t_{\cone(f)})^n}\lra\Phi(\cone(f))\ra
\cone(\Phi(t_{\cone(f)})^n)\ra \Phi(\cone(f))[1]$$ shows that for
large enough $n$ we get $\Phi(t_{\cone(f)})^n=0$ and so,
$\Phi(\cone(f))$ is a direct summand of
$\cone(\Phi(t_{\cone(f)})^n)$. But $\cone(\Phi(t_{\cone(f)})^n)$ is
in ${\mathcal C}_k^\circ$ and idempotents split in ${\mathcal
C}_k^\circ$. It follows that $Z$ is in ${\mathcal C}_k^\circ$.
\end{proof}

\subsection{Triangle degeneration implies categorical degeneration}

Following Keller~\cite{Kelleralgebraic}
a triangulated $k$-category $\mathcal C$ is algebraic if
$\mathcal C$ is the stable category of a Frobenius category.

\begin{Prop}\label{triangleimpliescat}
Let $\mathcal{C}_k^0$ be the category of compact objects of an
algebraic compactly generated triangulated $k$-category. If there is
a distinguished triangle $$Z\stackrel{\begin{pmatrix}v\\ u\end{pmatrix}}{\lra}Z\oplus
M\stackrel{\begin{pmatrix}h&j\end{pmatrix}}{\lra}N\lra Z[1] $$ in $\mathcal{C}_k^\circ$, with
a nilpotent endomorphism $v$ of $Z$,  then $M\leq_{cdeg}N$ (in the sense of
Definition~\ref{degendef}).
\end{Prop}

\begin{proof} Let $\mathcal{C}_k$ be a compactly generated algebraic triangulated
$k$-category such that $\mathcal{C}_k^\circ$ is the full subcategory
of its compact objects. By Keller~\cite[4.3 Theorem]{Kellerddc}, we
know that $\mathcal{C}_k$ is triangle-equivalent to the derived
category $\mathcal D (\mathcal A)$ of a small dg $k$-category
$\mathcal A$. We identify $\mathcal{C}_k=\mathcal D (\mathcal A)$ throughout
the rest of the proof.

Replacing each of the complexes $Z$, $M$ and $N$ by a homotopically
projective resolution in the homotopy category $\mathcal H(\mathcal
A)$ and adding suitable contractible complexes if necessary, we can
assume that the given triangle comes from a conflation
$$0\lra Z\stackrel{\begin{pmatrix}v\\ u\end{pmatrix}}{\lra}
Z\oplus M\stackrel{\begin{pmatrix}h&j\end{pmatrix}}{\lra} N\lra 0$$
in the complex category ${\mathcal C}({\mathcal A})$, where in addition $Z$,  $M$ and $N$
are finitely generated projective as graded $\mathcal A$-modules,
when we forget the differential.

Let us put $V:=k[[T]]$  and let $\mathcal A[[T]]$ be the dg
$V$-category with the same objects as $\mathcal A$, where
$$\text{Hom}_{\mathcal
A[[T]]} (A,A')=\text{Hom}_{\mathcal A}(A,A')[[T]]=\{\sum_{i=0}^\infty\alpha_iT^i\;:\;
\alpha_i\in\text{Hom}_{\mathcal A}(A,A')\text{, for all }i\geq 0\},$$
for all
$A,A'\in\mathcal A$, and the composition of morphisms is the obvious
one. We then denote by $\mathcal C (\mathcal A[[T]])$, $\mathcal H
(\mathcal{A}[[T]])$ and $\mathcal C_V:=\mathcal D(\mathcal A[[T]])$
their associated categories of complexes, homotopy and derived
category, respectively. We have an obvious functor
\begin{eqnarray*}?\hat{\otimes}
V:\mathcal C(\mathcal A)&\longrightarrow&\mathcal C (\mathcal
A[[T]])\\
M&\rightsquigarrow& M[[T]]
\end{eqnarray*}
where
$M[[T]]:(\mathcal A[[T]])^{op}\longrightarrow\mathcal C(V)$ acts as
$(M[[T]])(A)=M(A)[[T]]$ on objects and as
$M[[T]](\sum_{i=0}^\infty\alpha_iT^i )=\sum_{i=0}^\infty
M(\alpha_i)T^i$
on morphisms. The functor $?\hat{\otimes}V$ clearly
takes conflations to conflations and null-homotopic maps to
null-homotopic maps. We then get an induced triangulated functor
$$?\hat{\otimes}V:\mathcal H (\mathcal A)\longrightarrow\mathcal
H(\mathcal A[[T]]).$$ This latter functor clearly takes acyclic
complexes to acyclic complexes, so that we get a triangulated
functor $$\uar_k^V:\mathcal C_k=\mathcal D (\mathcal
A)\longrightarrow\mathcal D (\mathcal A[[T]])=:\mathcal C_V.$$ If
$N$ is any object of $\mathcal C_V$, then $N(A)$ is a graded
$V$-module, for each object $A$ of $\mathcal A$. Multiplication by
$T$ gives then a morphism $t_N(A):N(A)\longrightarrow N(A)$ of
graded $V$-modules of degree zero, for each object $A$ of $\mathcal A$. They
define a morphism $t_N:N\longrightarrow N$ in $\mathcal C_V$.
Moreover all the $t_N$ define a natural transformation
$t:\text{id}_{\mathcal C_V}\longrightarrow \text{id}_{\mathcal C_V}$
of triangulated functors. We claim that the data
$(\mathcal{C}_k,\mathcal C_V,\mathcal C_V^\circ ,\uar_k^V,t)$ give a
degeneration data for $\mathcal C_k^\circ =\mathcal D^c(\mathcal A)$
satisfying the requirements, where $\mathcal C_V^\circ =\mathcal D^c
(\mathcal A[[T]])$ is the category of compact objects of $\mathcal
C_V$. This will prove that $M\leq_{cdeg}N$.

There is an obvious `restriction of scalars' triangulated functor
$\phi :\mathcal D (\mathcal A[[T]])\longrightarrow\mathcal
D(\mathcal A)$. Moreover, if $X$ is an object of
$\mathcal C_k^\circ =\mathcal D^c(\mathcal A)$,
then $X$ is a direct summand of a finite iterated
extension of shifts of representable $\mathcal A$-modules $A^\wedge
=\text{Hom}_{\mathcal A}(?,A)$ (see \cite[Theorem 5.3]{Kellerddc}).
But we have an isomorphism in
$$(A^\wedge )\uar_k^V\cong\text{Hom}_{\mathcal A[[T]]} (?,A),$$
for each object $A$ of $\mathcal A$ (equivalently,  of $\mathcal
A[[T]]$). That is, the functor $\uar_k^V$ takes representable
$\mathcal A$-modules to representable $\mathcal A[[T]]$-modules, and
hence $(\mathcal C_k^\circ)\uar_k^V\subseteq\mathcal C_V^\circ$.

Let $M$ be any object of $\mathcal C_k$. The morphism
$$\phi (t_{M\uar_k^V}):\phi (M\uar_k^V)\longrightarrow \phi (M\uar_k^V),$$
when evaluated at an object $A$ of $\mathcal A$ (or $\mathcal
A[[T]]$), gives the morphism of $k$-modules $$M(A)[[T]]\longrightarrow
M(A)[[T]]$$ given by multiplication by $T$. We get a retraction
$$\pi_M(A):M(A)[[T]]\longrightarrow M(A)[[T]]$$
for this map
which takes
$\sum_{i=0}^\infty m_iT^i$
to $\sum_{i=1}^\infty m_iT^{i-1}$.
It is routine to see that this gives a
retraction $\pi_M:\phi (M\uar_k^V)\longrightarrow \phi
(M\uar_k^V)$ for $\phi (t_{M\uar_k^V})$. Moreover, in the abelian category
$\mathcal C(\mathcal A[[T]])$, we have an exact sequence
$$0\rightarrow M[[T]]\stackrel{\cdot T}{\longrightarrow} M[[T]]\longrightarrow
M[[T]]/TM[[T]]\rightarrow 0.$$
We then get a distinguished triangle
$$M\uar_k^V\stackrel{t_{M\uar_k^V}}{\longrightarrow} M\uar_k^V
\longrightarrow M[[T]]/TM[[T]]\longrightarrow M\uar_k^V [1]$$
 in $\mathcal D(\mathcal A[[T]])$.
It is clear that
$$\text{cone}(\phi (t_{M\uar_k^V}))\cong\phi(M[[T]]/TM[[T]])$$
is isomorphic to $M$ in $\mathcal D(\mathcal{A})$. Then all
conditions of Definition \ref{degendatadef} are satisfied.

We now check that $M\leq_{cdeg}N$ with respect to the above-defined
degeneration data. We first claim that
$\begin{pmatrix}v+T\\ u\end{pmatrix}:Z[[T]]\longrightarrow Z[[T]]\oplus M[[T]]$ is an inflation
(=admissible monomorphism) in the exact category $\mathcal
C(\mathcal A[[T]])$. For this we need to get a retraction for this
map as in the category of graded $\mathcal A[[T]]$-modules. We look
for a retraction of the form
$\begin{pmatrix}\left(\sum_{i=0}^\infty h_iT^i\right)&\left(\sum_{i=0}^\infty j_iT^i\right)\end{pmatrix}$,
where $h_i:Z\longrightarrow Z$
and $j_i:M\longrightarrow Z$ are morphisms of graded $\mathcal
A$-modules, for all $i\geq 0$. We already know that $\begin{pmatrix}v\\ u\end{pmatrix}$
allows such a retraction $\begin{pmatrix}h&j\end{pmatrix}$ in the category of graded $\mathcal
A$-modules. The equation
$$\begin{pmatrix}\left(\sum_{i=0}^\infty h_iT^i\right)&\left(\sum_{i=0}^\infty j_iT^i\right)\end{pmatrix}
\cdot\begin{pmatrix}{v+T}\\ {u}\end{pmatrix}=\text{id}_{Z[[T]]}$$
translates into
\begin{eqnarray*}
j_0u+h_0v&=&\text{id}_Z\\
h_0+j_1u+h_1v&=&0\\
h_1+j_2u+h_2v&=&0\\
\dots&\dots&\dots\\
h_{i-1}+j_iu+h_iv&=&0\\
\dots&\dots&\dots
\end{eqnarray*}
and so $h_0=h$ and $j_0=j$ solves the first equation.
Since $$\text{id}_Z=\begin{pmatrix}h&j\end{pmatrix}\cdot\begin{pmatrix}v\\ u\end{pmatrix},$$
we get
$$-h=\begin{pmatrix}-h^2&-hj\end{pmatrix}\cdot \begin{pmatrix}v\\ u\end{pmatrix}$$
and if we know $h_{i-1}$, then
$$-h_{i-1}=\begin{pmatrix}-h_{i-1}h&-h_{i-1}j\end{pmatrix}\cdot \begin{pmatrix}v\\ u\end{pmatrix},$$
so that we may put $h_i:=-h_{i-1}h$ and $j_i:=-h_{i-1}j$. This gives
$h_i=(-1)^ih^{i+1}$ and $j_i=(-1)^ih^{i}j$ by an obvious induction. We
then obtain a conflation in $\mathcal C(\mathcal A[[T]])$
$$(\dagger)\;\;\;0\lra Z[[T]]
\stackrel{\begin{pmatrix}{v+T}\\ {u}\end{pmatrix}}{\lra} Z[[T]]\oplus M[[T]]\stackrel{}{\lra} Q\lra 0$$
and we get a commutative diagram
$$\begin{array}{ccccccc}
0\ra& Z[[T]]&\stackrel{\begin{pmatrix}{v+T}\\ {u}\end{pmatrix}}{\ra}&
Z[[T]]\oplus M[[T]]&\stackrel{}{\ra}& Q&\ra 0\\
&\phantom{\mbox{\scriptsize$T$}}\dar \mbox{\scriptsize$T$}&&
\phantom{\mbox{\scriptsize $\left(\begin{array}{cc}T&0\\0&T\end{array}\right)$}}
\dar\mbox{\scriptsize $\left(\begin{array}{cc}T&0\\0&T\end{array}\right)$}
&&\phantom{\mbox{\scriptsize$T$}}\dar\mbox{\scriptsize$T$}\\
0\ra& Z[[T]]&\stackrel{\begin{pmatrix}{v+T}\\ {u}\end{pmatrix}}{\ra}& Z[[T]]\oplus M[[T]]&
\stackrel{}{\ra}& Q&\ra 0
\end{array}$$

The canonical functor
$p:\mathcal{C}_V^\circ\longrightarrow\mathcal{C}_V^\circ [t^{-1}]$
maps the endomorphism $v+T$ of $Z[[T]]$ onto an invertible
endomorphism of $p(Z[[T]])$ since $v$ is nilpotent and
multiplication by $T$ is an isomorphism in $\mathcal{C}_V^\circ
[t^{-1}]$. The triangle given by the conflation ($\dagger$) above is then mapped by $p$ onto
a split triangle in $\mathcal{C}_V^\circ[t^{-1}]$ such that the
first component of the morphism
$$p(Z[[T]])
\stackrel{\begin{pmatrix} p(v+T)\\p(u)\end{pmatrix}}
{\longrightarrow}
p(Z[[T]])\oplus p(Z[[T]])$$ is an
isomorphism. We then get that $p(M\uar_k^V)=p(M[[T]])$ is isomorphic
to $p(Q)$ in $\mathcal{C}_V^\circ [t^{-1}]$.

\medskip

Moreover, multiplication by $T$ gives a monomorphism $X[[T]]\ra
X[[T]]$  in the (abelian) category of complexes ${\mathcal
C}({\mathcal A}[[T]])$ of $\mathcal A[[T]]$-modules, for each
$X\in\mathcal C(\mathcal A)$. In addition,  the forgetful functor ${\mathcal
C}({\mathcal A}[[T]])\longrightarrow\text{Gr}-{\mathcal A}[[T]]$ is exact, where
$\text{Gr}-{\mathcal A}[[T]]$ denotes the category of graded $\mathcal{A}[[T]]$-modules.
This together with the split monomorphic condition of $v+T$ in $\text{Gr}-\mathcal{A}[[T]]$
imply that the map $Q\stackrel{T}{\longrightarrow}Q$ in the above diagram is also a monomorphism
in ${\mathcal C}({\mathcal A}[[T]])$. Applying now the kernel-cokernel lemma to that diagram,
we get a short exact sequence
$$0\lra Z[[T]]/TZ[[T]]\lra Z[[T]]/TZ[[T]]\oplus M[[T]]/TM[[T]]\lra Q/TQ\lra 0$$
in ${\mathcal C}({\mathcal A}[[T]])$. Identifying $X[[T]]/TX[[T]]$ with $X$,
for each object $X$ of $\mathcal C(\mathcal A)$,  the first map is again
$\begin{pmatrix}v\\ u\end{pmatrix}$, and hence $Q/TQ\simeq N$ in ${\mathcal C}({\mathcal
A})$. Hence $\Phi(\cone(t_Q))\simeq N$ in  ${\mathcal D}^c({\mathcal
A})={\mathcal C}_k^\circ$. This shows $M\leq_{cdeg} N$.
\end{proof}

\subsection{The main result}

Our main result is the following.

\begin{Theorem}\label{main} Let $k$ be a commutative ring and
let ${\mathcal C}_k^\circ$ be a triangulated $k$-category with split
idempotents. If $M$ and $N$ are objects of $\mathcal{C}_k^\circ$ and
$M$ degenerates to $N$ in the categorical sense, then there is a
distinguished triangle
$$Z\stackrel{\begin{pmatrix}v\\ u\end{pmatrix}}{\lra} Z\oplus M\lra N\lra Z[1]$$
with nilpotent endomorphism $v$. When $\mathcal{C}_k^\circ$ is
equivalent to the category of compact objects of a compactly
generated algebraic triangulated $k$-category, the converse is also
true.
\end{Theorem}

\begin{proof}
This is an immediate consequence of Proposition~\ref{catimpliestriangle}
and Proposition~\ref{triangleimpliescat}.
\end{proof}

\begin{Rem}
For any abelian category $\mathcal A$ define
$M\leq_{ses} N$ if there is a short exact sequence
$0\ra Z\stackrel{\begin{pmatrix}v\\ u\end{pmatrix}}{\ra}Z\oplus M\ra N\ra 0$
without any further hypothesis on $v$.
If $\mathcal A$ admits countable direct sums we showed in \cite[Example 2.2.(1)]{JSZpartord},
then $M\leq_{ses} N$ and $N\leq_{ses} M$ for
each pair of objects $N$ and $M$ of $\mathcal A$. The example transposes immediately to the
triangulated situation.
We insist in the fact that there the endomorphism $v$
is not nilpotent, as we supposed in the hypotheses of Theorem~\ref{main}.
If the endomorphism algebras of each of the concerned objects
is artinian, then the condition on $v$ to be nilpotent is superfluous in Theorem~\ref{main}.
\end{Rem}

\begin{Rem}
The definition that an object $M$ degenerates to an object $N$ of a triangulated category ${\mathcal D}^\circ_k$ in the categorical sense depends a priori on the category ${\mathcal D}^\circ_k$. Do we get the same concept if ${\mathcal D}^\circ_k$ is a full subcategory of a triangulated category ${\mathcal C}^\circ_k$? This question was answered at least partially in \cite[Proposition 11]{Osakaproceedings}.

Let $k$ be a field and let ${\mathcal C}_k^\circ$ be the category of compact objects in an algebraic compactly generated triangulated $k$-category. If ${\mathcal D}^\circ_k$ is a full triangulated subcategory of ${\mathcal C}_k^\circ$, then for all objects
$M$ and $N$ of ${\mathcal D}_k^\circ$ we get that $M\leq_{cdeg}N$ with respect to
${\mathcal D}_k^\circ$ if and only if $M\leq_{cdeg}N$ in ${\mathcal C}_k^\circ$.
\end{Rem}

\section{Some concluding remarks on triangulated degeneration}

\label{conclusionsection}

\begin{Def}
Let $\mathcal{T}$ be any triangulated category. We shall denote by
$\preceq_{cdeg}$ (resp. $\preceq_\Delta$, resp. $\preceq_{\Delta*}$) the smallest transitive
relation in the class of objects of $\mathcal{T}$ satisfying that if
$X\leq_{cdeg}Y$ (resp. there exists a distinguished triangle
$Z\stackrel{\begin{pmatrix} v\\ u
\end{pmatrix}}{\longrightarrow}Z\oplus X\longrightarrow Y\longrightarrow Z[1]$,
resp. there is such a triangle,  with $v$ a
nilpotent endomorphism of $Z$), then $X\preceq_{cdeg}Y$
(resp. $X\preceq_\Delta Y$, resp. $X\preceq_{\Delta*}Y$).
\end{Def}

\begin{Rem} \label{rem.Jensen-Su-Zimmermann}
We recall that Jensen, Su and the second author showed in \cite{JSZpartord}
that if $\mathcal T$ is a triangulated category in which idempotent endomorphisms
split, and in which the endomorphism ring of each object is artinian,
then $\preceq_\Delta$ and $\preceq_{\Delta*}$ coincide and they can be
defined in one step, i.e., $X\preceq_\Delta Y$ if, and only if, there is
a distinguished triangle
$$Z\stackrel{
\begin{pmatrix} v\\ u\end{pmatrix}}{\longrightarrow}
Z\oplus X\longrightarrow Y\longrightarrow Z[1].$$
If moreover $\mathcal T$ is a
$k$-linear category over a commutative ring $k$ such that the
composition length (as $k$-module)
of the homomorphism space between two objects is always finite, and if for
each two objects $X$ and $Y$ in $\mathcal T$ there is $n\in\Z$ so that
$\text{Hom}_{\mathcal T}(X,Y[n])=0$, then $\preceq_\Delta$ is also
antisymmetric. In this case,  when $\mathcal T$ is skeletally small,
the relation $\preceq_\Delta =\preceq_{\Delta*}$ is a partial order in
the set of isomorphism classes of objects of $\mathcal{T}$.
\end{Rem}

\begin{Prop} Let $k$ be a commutative ring.
Let ${\mathcal T}$ be a triangulated category with split idempotents and suppose that
for all objects $X$ and $Y$ the $k$-module
$\text{Hom}_{\mathcal T}(X,Y)$ is of finite composition length.
If $\text{Hom}_{{\mathcal T}}(M,M[1])=0$ then
$M$ is a minimal object
with respect to the relation $\preceq_{\Delta}=\preceq_{\Delta*}$. More precisely, under these conditions
$N\preceq_{\Delta}M$ implies $N\simeq M$.
\end{Prop}

\begin{proof}
Note that for each object $X$ of $\mathcal C$ we have that
$\text{End}_{\mathcal T}(X)$ is an Artinian $k$-algebra. Indeed
$\text{End}_\mathcal{T}(X)$ is a
$k$-module of finite length and then
$I:=\text{ann}_k(\text{End}_\mathcal{T}(X))$ is an ideal of $k$ such
that $k/I$ has finite length as $k$-module, and so $k/I$ is Artinian.
It follows that $\text{End}_\mathcal{T}(X)$
is an Artinian $k$-algebra since it has finite length as a module over
the artinian commutative ring $k/I$.

Let
$$Z\stackrel{\begin{pmatrix} u\\ v \end{pmatrix}}{\lra} N\oplus Z\ra M\ra Z[1]$$
be a distinguished triangle, where we assume without loss of generality that
$v$ is nilpotent (see Remark \ref{rem.Jensen-Su-Zimmermann}).
We need to show that $M\simeq N$.
Apply $\text{Hom}_{{\mathcal T}}(M,-)=:(M,-)$
to the distinguished triangle and part of the long exact sequence becomes
$$(M,Z)\ra (M,N)\oplus (M,Z)\ra (M,M)\ra (M,Z[1])\ra (M,N[1])\oplus (M,Z[1])\ra 0
$$
using that $(M,M[ 1])=0.$
This shows that $(M,N[ 1])=0$ and that the sequence
$$(M,Z)\longrightarrow (M,N)\oplus (M,Z)\longrightarrow(M,M)\rightarrow 0$$
is exact. Then the given distinguished triangle splits since $\text{id}_M$ is in the image of the
map $(M,N)\oplus (M,Z)\longrightarrow(M,M)$.
This shows that $N\oplus Z\simeq M\oplus Z$.
By \cite[Theorem A.1]{XiaoWuChen} we obtain that $\mathcal T$ is a Krull-Schmidt category, and hence
$N\cong M$.
\end{proof}

\begin{Rem}
If $\mathcal T$ is a triangulated category and if
$Z\ra Z\oplus M\ra N\ra Z[1]$ is a split distinguished triangle such that the induced endomorphism $v$
on $Z$ is in the Jacobson radical of the endomorphism ring of $Z$, then we can show that
$M$ is isomorphic to $N$, even if $\mathcal T$ is not a Krull-Schmidt category.
\end{Rem}

\begin{Rem}
As observed by Yoshino~\cite[Remark 4.6]{Yoshinostable}
there cannot be a maximal element with respect to $\leq_\Delta$. Indeed,
$$X\lra 0\lra X[1]\lra X[1]$$ is a distinguished triangle, and hence
also for any object $U$ of $\mathcal T$ we get that
$$X\stackrel{\begin{pmatrix}0\\ 0\end{pmatrix}}{\lra} U\oplus X\lra U\oplus X\oplus X[1]\lra X[1]$$
is a distinguished triangle.
Choosing $Z:=X$ this shows $U\leq_\Delta (U\oplus X\oplus X[1]),$
i.e.  that the object $U$ degenerates to
$U\oplus X\oplus X[1]$ for all objects $X$.
Observe that here the endomorphism on $X$ on the left of the distinguished triangle
is actually $0$, hence nilpotent, so that we are in the situation of Theorem~\ref{main}.
\end{Rem}

For a commutative Noetherian ring $R$, we will denote
by $R-fl$ the category of finite length $R$-modules.

\begin{Lemma} \label{lem.right exact are exact}
Let $R$ be a commutative ring and let  $U\stackrel{\begin{pmatrix}\alpha \\
\beta
\end{pmatrix}}{\longrightarrow}U\oplus V\stackrel{\begin{pmatrix} \gamma &
\delta \end{pmatrix}}{\longrightarrow}W\rightarrow
0$ be an exact sequence of finitely generated $R$-modules.
Suppose that either one of the following conditions holds:
\begin{enumerate}
\item $R$ is Noetherian.
\item The $R$-modules $U$, $V$ and $W$ have finite length.
\end{enumerate}
If $\text{length}(\text{Hom}_R(V,Y))=\text{length}(\text{Hom}_R(W,Y))$,
for all $R$-modules of finite length $Y$, then $\begin{pmatrix}\alpha \\
\beta\end{pmatrix}$
is a monomorphism and the short exact sequence
$0\rightarrow U\stackrel{\begin{pmatrix}\alpha \\ \beta
\end{pmatrix}}{\longrightarrow}U\oplus V\stackrel{\begin{pmatrix} \gamma &
\delta \end{pmatrix}}{\longrightarrow}W\rightarrow
0$ remains exact when we apply the functor $\text{Hom}_R(?,Y)$, for
any $R$-module $Y$ of finite length.
Moreover, under condition \mbox{\rm (2)} the sequence is split.
\end{Lemma}

\begin{proof}
Let us fix an $R$-module $Y$ of finite length. When we apply the contravariant functor
$\text{Hom}_R(?,Y)$, we get an exact sequence in $R-fl$

$$0\rightarrow\text{Hom}_R(W,Y)\stackrel{\begin{pmatrix} \gamma^*\\ \delta^*\end{pmatrix}}{\longrightarrow}\text{Hom}_R(U,Y)\oplus\text{Hom}_R(V,Y)
\stackrel{\begin{pmatrix} \alpha^* & \beta^*
\end{pmatrix}}{\longrightarrow}\text{Hom}_R(U,Y).$$
By comparison of composition lengths,
we get that
$$\text{length}(\text{Hom}_R(U,Y))+\text{length}(\text{Hom}_R(V,Y))\leq
\text{length}(\text{Hom}_R(U,Y))+
\text{length}(\text{Hom}_R(W,Y))$$
and hence $$\text{length}(\text{Hom}_R(V,Y))\leq \text{length}(\text{Hom}_R(W,Y))$$
with equality  if and only if
$\begin{pmatrix}
\alpha^* & \beta^*
\end{pmatrix}$ is surjective.

Putting $f=\begin{pmatrix}\alpha \\ \beta
\end{pmatrix}$, the proof reduces to check
that if $f:M\longrightarrow N$ is a morphism of finitely generated
$R$-modules such that
$f^*:\text{Hom}_R(N,Y)\longrightarrow\text{Hom}_R(M,Y)$ is an
epimorphism for all $Y\in R-fl$, then $\ker(f)=0$.

Suppose that $U$, $V$ and $W$ are modules of finite length. Then
we can take $Y=M$ and we get that $\textrm{id}_M$ is in the image of
$f^*$, which implies that $f$ is a split monomorphism. This implies that
$\ker(f)=0$.

Suppose now that $R$ is Noetherian.
By localizing at maximal ideals, we can assume without
loss of generality that $R$ is a local ring. In such case, if
$\mathbf{m}$ is the maximal ideal and we take $Y=M/\mathbf{m}^nM$,
the fact that the projection $p:M\longrightarrow M/\mathbf{m}^nM$ is
in the image of $f^*$ implies that
$\ker(f)\subseteq\mathbf{m}^nM$. It follows that
$\ker(f)\subseteq\bigcap_{n\geq 0}\mathbf{m}^nM$, so that
$\ker(f)=0$ (see \cite[Theorem 8.10]{Matsumura}).
\end{proof}

\begin{Lemma} \label{lem.crucial for equivalence relation}
Let $k$ be a commutative ring, let $\mathcal{T}$ be a  triangulated
$k$-category, let $H:\mathcal{T}\longrightarrow k-\text{mod}$ be a
cohomological functor, where $k-\text{mod}$ denotes the category of
finitely generated $k$-modules. Put $H^j=H\circ (?[j])$, for each
integer $j$, and suppose that either one of the following conditions
holds:
\begin{enumerate}
\item $k$ is Noetherian and, for each object $X$ of $\mathcal{T}$, the $k$-module
$H(X)$ is  finitely generated
and there is an integer $n=n_X$ such that $H^j(X)=0$ for all $j>n$;
\item For each object $X$ of $\mathcal{T}$, the $k$-module $H(X)$ is
of finite length and one of the following subconditions holds:
\begin{enumerate}
\item For each object $X$ of $\mathcal{T}$, there is an integer
$n=n_X$ such that $H^j(X)=0$, for all $j>n$;
\item For each object $X$ of $\mathcal{T}$, there is an integer
$n=n_X$ such that $H^j(X)=0$, for all $j<n$
\end{enumerate}
\end{enumerate}
If $Z\stackrel{\begin{pmatrix} v \\ u
\end{pmatrix}}{\longrightarrow}Z\oplus M\longrightarrow N\longrightarrow
Z[1]$ is a distinguished triangle in $\mathcal{T}$ and $N\preceq_\Delta M$,
then the long exact sequence
associated to the cohomological functor $H$ gives short exact sequences
$$0\rightarrow H^j(Z)\longrightarrow H^j(Z)\oplus
H^j(M)\longrightarrow H^j(N)\rightarrow 0,$$ for all integers $j$,
which remain exact when we apply the functor $\text{Hom}_k(?,Y)$,
for any $k$-module $Y$ of finite length. In particular, they are
all split exact under condition \mbox{\rm (2)}.
\end{Lemma}

\begin{proof}
By definition of $\preceq_\Delta$, we have a sequence
$M_0,M_1,\dots,M_r$ of objects of $\mathcal{T}$, with $M_0=M_r=M$ and
$M_1=N$, together with distinguished triangles
$$Z_i\longrightarrow Z_i\oplus M_i\longrightarrow M_{i+1}\longrightarrow Z_i[1]$$
($i=0,1,\dots,r-1$).

We first assume that condition (1) or condition (2).(a) hold. We can take
$n\in\mathbb{Z}$ such that $H^j(Z_i)=0=H^j(M_i)$, for all
$i=0,1,\dots,r-1$ and all $j>n$. We then get exact sequences of
finitely generated $k$-modules
$$H^n(Z_i)\longrightarrow H^n(Z_i)\oplus H^n(M_i)\longrightarrow H^n(M_{i+1})\rightarrow 0$$
($i=0,1,\dots,r-1$). If $Y$ is any $k$-module of finite length and we
apply the functor $\text{Hom}_k(?,Y)$ then, counting composition
lengths as in the proof of Lemma~\ref{lem.right exact are exact}, we get that
$$\text{length}(\text{Hom}_k(H^n(M_i),Y))\leq\text{length}(\text{Hom}_k(H^n(M_{i+1}),Y))$$
for $i=0,1,\dots,r-1$. We then get that these inequalities are
equalities since
$$\text{length}(\text{Hom}_k(H^n(M_0),Y))=\text{length}(\text{Hom}_k(H^n(M_r),Y)).$$
Lemma \ref{lem.right exact are exact} implies that we get
short exact sequences
$$0\rightarrow H^n(Z_i)\longrightarrow H^n(Z_i)\oplus H^n(M_i)
\longrightarrow H^n(M_{i+1})\rightarrow 0,$$
which remain exact when we apply $\text{Hom}_k(?,Y)$ for any
$k$-module $Y$ of finite length, and which split under hypothesis (2).(a).

The long exact sequence associated to the cohomological functor $H$ then gives exact sequences
$$H^{n-1}(Z_i)\longrightarrow H^{n-1}(Z_i)\oplus H^{n-1}(M_i)
\longrightarrow H^{n-1}(M_{i+1})\rightarrow 0$$
($i=0,1,\dots,r-1$). Repeating now iteratively  for all $j<n$ the
argument used for $n$, we obtain the result.

We now assume that condition (2).(b) holds and fix an integer $n$ such
that $H^j(Z_i)=0=H^j(M_i)$, for all $j<n$ and $i=0,1,\dots,r-1$. We then get  exact
sequences
$$0\rightarrow H^n(Z_i)\longrightarrow H^n(Z_i)\oplus H^n(M_i)\longrightarrow H^n(M_{i+1})$$
($i=0,1,\dots,r-1$). It follows that
$\text{length}(H^n(M_i))\leq\text{length}(H^n(M_{i+1}))$, for each
$i=0,1,\dots,r-1$, and these inequalities are then equalities since
$H^n(M_0)=H^n(M_r)$. We  get  short exact sequences
$$(*)\;\;\;0\rightarrow H^n(Z_i)\longrightarrow H^n(Z_i)\oplus H^n(M_i)\longrightarrow H^n(M_{i+1})\rightarrow
0. $$
When we take a $k$-module of finite length $Y$, apply
the functor $\text{Hom}_k(?,Y)$ to these sequences and count
composition lengths to obtain that
$$\text{length}(\text{Hom}_k(H^n(M_i),Y))\leq
\text{length}(\text{Hom}_k(H^n(M_{i+1}),Y)),$$
for each
$i=0,1,\dots,r-1$. Again, these inequalities are equalities and, using
Lemma \ref{lem.right exact are exact}, we conclude that all the
exact sequences $(*)$ split. Finally, using an argument dual to the
one followed under conditions (1) or (2).(a),  one easily shows by
induction on $s\geq 0$ that the long exact sequence associated to
$H$ gives split short exact sequences
$$0\rightarrow H^{n+s}(Z_i)\longrightarrow H^{n+s}(Z_i)\oplus
H^{n+s}(M_i)\longrightarrow H^{n+s}(M_{i+1})\rightarrow 0,$$ for
$i=0,1,\dots,r-1$.
\end{proof}

We are now ready to prove the main result of this section.

\begin{Theorem} \label{theor.triangle-degen. partial order}
Let $k$ be a commutative ring and let $\mathcal{T}$ be a skeletally
small triangulated $k$-category with split idempotents. The
following assertions hold:

\begin{enumerate}
\item Suppose that $\text{Hom}_\mathcal{T}(X,Y)$ is a $k$-module of finite length,
for all objects $X,Y$ of $\mathcal{T}$, and either one of the
following two conditions hold for such objects:

\begin{enumerate}
\item There is an integer $n=n_{XY}$ such that
$\text{Hom}_\mathcal{T}(X,Y[j])=0$, for all $j>n$;
\item There is an integer $n=n_{XY}$ such that
$\text{Hom}_\mathcal{T}(X,Y[j])=0$, for all $j<n$.
\end{enumerate}
Then $\preceq_\Delta =\preceq_{\Delta*}$ is a partial order in the set of isomorphism classes of
objects of $\mathcal{T}$.

\item Suppose that $k$ is Noetherian, that
$\text{Hom}_\mathcal{T}(X,Y)$ is a finitely generated $k$-module,
for all objects $X,Y$ of $\mathcal{T}$, and that there is an $n=n_{XY}$ such that
$\text{Hom}_{\mathcal T}(X,Y[j])=0$ for all $j>n$.
If $M\preceq_\Delta N$ and $N\preceq_\Delta M$
then there is an object $Z$ of $\mathcal{T}$ such that $M\oplus Z$
and $N\oplus Z$ are isomorphic.
\end{enumerate}
\end{Theorem}

\begin{proof}
For the equality $\preceq_\Delta =\preceq_{\Delta*}$ in Assertion (1),
see Remark \ref{rem.Jensen-Su-Zimmermann}.
Both in assertion (1) and (2), the relation $\preceq_\Delta$ is reflexive and transitive.
Let us assume that $M\preceq_\Delta N$ and $N\preceq_\Delta M$. We
then have a sequence of objects $M_0,M_1,\dots,M_r$ in $\mathcal{T}$
such that $M_0=M=M_r$ and $M_s=N$, for some $s=0,1,\dots,r$ and we
have distinguished triangles
$$(**)\;\;\;\;\;Z_i\longrightarrow Z_i\oplus M_i\longrightarrow M_{i+1}\longrightarrow Z_i[1] $$
($i=0,1,\dots,r-1$). Let $X$ be any object of $\mathcal{T}$ and
consider the cohomological functor
$\text{Hom}_\mathcal{T}(X,?):\mathcal{T}\longrightarrow
k-\text{mod}$. By Lemma~\ref{lem.crucial for equivalence relation},
under any of the assumptions, we have exact sequences
$$0\rightarrow\text{Hom}_\mathcal{T}(X,Z_i[j])\longrightarrow
\text{Hom}_\mathcal{T}(X,Z_i[j])\oplus\text{Hom}_\mathcal{T}(X,M_i[j])
\longrightarrow\text{Hom}_\mathcal{T}(X,M_{i+1}[j])\rightarrow 0,$$
for all $j\in\mathbb{Z}$ and $i=0,1,\dots,r-1$. Fixing $i$ and putting
$X=M_{i+1}$ and $j=0$, we get that all the triangles $(**)$ split.
Therefore $M_i\oplus Z_i\cong M_{i+1}\oplus Z_i$, for all
$i=0,1,\dots,r-1$, which easily implies  for the
object $Z=Z_0\oplus\dots\oplus Z_{s-1}$ of $\mathcal{T}$
that $M\oplus Z\cong N\oplus Z$.

When $\text{Hom}_\mathcal{T}(X,Y)$ is a $k$-module of finite length,
for all objects $X,Y$ of $\mathcal{T}$, we have that $\mathcal{T}$
is a Krull-Schmidt category. Indeed $\text{End}_\mathcal{T}(X)$ is a
$k$-module of finite length and then
$I:=\text{ann}_k(\text{End}_\mathcal{T}(X))$ is an ideal of $k$ such
that $k/I$ is Artinian. It follows that $\text{End}_\mathcal{T}(X)$
is an Artinian $k$-algebra, for each object $X$ of $\mathcal{T}$. By
\cite[Theorem A.1]{XiaoWuChen} we get that $\mathcal T$ is a Krull-Schmidt category.
Hence, from the isomorphism $M\oplus Z\cong N\oplus Z$
we get that $M\cong N$.
\end{proof}

\begin{Example}
We recall Example 2.2.(2) of \cite{JSZpartord}.
Let $G$ be the generalised quaternion group of order $32$.
Swan gives an ideal $\mathfrak a$ of $\Z G$ such that
$\Z G\oplus {\mathfrak a}\simeq \Z G\oplus\Z G$ as $\Z G$-modules, and such that
${\mathfrak a}\not\simeq \Z G$. This implies that
$\Z G\leq_\Delta{\mathfrak a}\leq_\Delta \Z G$ when we choose ${\mathcal T}=D^b(\Z G)$.
\end{Example}

Recall that a ring $A$ is called a an \emph{Artin algebra} (resp. a
\emph{Noether algebra}) when its center is Artinian (resp.
Noetherian) and $A$ is finitely generated as a module over that
center.

\begin{Cor}
If $A$ is an Artin algebra and $\mathcal{T}=D^b(A)$, then the
relations $\preceq_{cdeg}$,  $\preceq_\Delta$ and $\preceq_{\Delta*}$ coincide in $D^b(A)$
and are partial orders in the set of isomorphism classes of objects of
$D^b(A)$.
\end{Cor}

\begin{proof}
By Remark \ref{rem.Jensen-Su-Zimmermann}, we know that $\preceq_\Delta =\preceq_{\Delta*}$.
Let $M,N$ be any objects of $D^b(A)$ and suppose that
$M\preceq_{cdeg}N$. By definition of $\preceq_{cdeg}$, there is a
sequence of objects $M_0,M_1,\dots,M_r$ in $D^b(A)$ such that $M_0=M$
and  $M_r=N$ and such that $M_i\leq_{cdeg}M_{i+1}$ for all $i=0,1,\dots,r-1$.
By Theorem \ref{main}, we get that
$M\preceq_\Delta N$. For the converse implication, note that
$D^b(A)$ is equivalent to the subcategory of compact objects of a
compactly generated algebraic triangulated category (see
\cite{Krauseschemes}). Then, again by Theorem \ref{main}, we get
that  $M\preceq_\Delta N$ implies $M\preceq_{cdeg}N$.

Finally, $\mathcal{T}=D^b(A)$ satisfies condition (1).(b) of
Theorem~\ref{theor.triangle-degen. partial order} and, hence,
$\preceq_\Delta =\preceq_{cdeg}$ is a partial order in the set of
isomorphism classes of objects of $D^b(A)$.
\end{proof}

Recall that a dg algebra $A$ is \emph{homologically upper  bounded}
when there is an integer $m$ such that $H^n(A)=0$, for all $n>m$.

\begin{Cor}
Let  $k$ be Noetherian commutative and let
$A=\oplus_{n\in\mathbb{Z}}A^n$ be a dg $k$-algebra such that
$H^n(A)$ is a finitely generated $k$-module, for all
$n\in\mathbb{Z}$. Suppose that $A$ is homologically upper bounded
(e.g., a Noether algebra viewed as a dg algebra concentrated in
degree 0) and let $\textrm{per}(A)$ be the category of compact objects
of $D(A)$. The relations $\preceq_{\Delta^*}$ and $\preceq_{cdeg}$
coincide in $\textrm{per}(A)$. Moreover, if $M,N$ are objects of $\textrm{per}(A)$
such that $M\preceq_\Delta N$ and $N\preceq_\Delta M$ (in particular,
if $M\preceq_{cdeg}N$ and $N\preceq_{cdeg}M$), then the following
assertions hold:
\begin{enumerate}
\item There is $Z\in\textrm{per}(A)$ such that $Z\oplus M$ and $Z\oplus
N$ are isomorphic in $D(A)$;
\item The $k_\wp$-modules $H^j(M)_\wp$ and
$H^j(N)_\wp$ are isomorphic, for all $j\in\mathbb{Z}$ and all
prime ideals $\wp$ of $k$.
\end{enumerate}
\end{Cor}

\begin{proof}
That $\preceq_{\Delta^*}$ and $\preceq_{cdeg}$ coincide in
$\text{per}(A)$ is a direct consequence of Theorem \ref{main} and the definition of these relations.
Assertion (1)
will  follow from Theorem \ref{theor.triangle-degen. partial order}
once we check that if $X,Y$ are objects of $\textrm{per}(A)$, then
$\textrm{Hom}_{D(A)}(X,Y)$ is a finitely generated $k$-module and
$\textrm{Hom}_{D(A)}(X,Y[j])=0$ for $j>>0$. Indeed, each object of
$\textrm{per}(A)$ is a direct summand of a finite iterated extension
of shifts $A[n]$, with $n\in\mathbb{Z}$ (see \cite[Theorem
5.3]{Kellerddc}). This reduces the proof to the case when $X=A[m]$
and $Y=A[n]$, for some $m,n\in\mathbb{Z}$. We have that
$\textrm{Hom}_{D(A)}(A[m],A[n])\cong H^{n-m}(A)$, which is a finitely
generated $k$-module. On the other hand,  we have an equality
\begin{eqnarray*}
\textrm{Hom}_{D(A)}(A[m],A[n][j])&=&\textrm{Hom}_{D(A)}(A[m],A[n+j])\\
&\cong&\textrm{Hom}_{D(A)}(A,A[n-m+j])\\
&\cong& H^{n-m+j}(A).
\end{eqnarray*}
This is zero for $j>>0$ due to the homologically upper bounded
condition of $A$.

Given an isomorphism   $Z\oplus M\cong Z\oplus N$, we get that
$\text{length}(\text{Hom}_k(H^j(M),Y))=\text{length}(\text{Hom}_k(H^j(M),Y))$,
for each $k$-module $Y$ of finite length. Assertion (2) then follows
from \cite[Theorem 2.2]{YoshinoCM}.
\end{proof}


\begin{thebibliography}{88}


\bibitem{XiaoWuChen}
Xiao-Wu Chen, Yu Ye, and Pu Zhang, {\em Algebras of Derived Dimension Zero},
Communications in Algebra {\bf 36} (2008) 1-10.

\bibitem{GabrielZisman}
Pierre Gabriel and Michel Zisman, {\sc Calculus of Fractions and Homotopy Theory},
Springer Verlag Heidelberg 1967.


\bibitem{HuisgenSaorin}
Birge Huisgen-Zimmermann and Manuel Saorin,
{\em Geometry of Chain Complexes and Outer Automorphisms under Derived Equivalences},
Transactions of the American Mathematical Society
{\bf 353} (2001) 4757-4777.

\bibitem{JMS}
Bernt Tore Jensen, Dag Madsen and Xiuping Su,
{\em Degeneration of A-infinity modules },
Transactions of the American Mathematical Society {\bf 361}
(2009) 4125-4142.

\bibitem{JSZdegen}
Bernt Tore Jensen, Xiuping Su and Alexander Zimmermann,
{\em Degenerations for derived categories},
Journal of Pure and Applied Algebra {\bf 198} (2005) 281-295.

\bibitem{JSZpartord}
Bernt Tore Jensen, Xiuping Su and Alexander Zimmermann, {\em
Degeneration-like orders for triangulated categories},
Journal of Algebra and its Applications {\bf 4} (2005) 587-597.

\bibitem{Kellerddc}
Bernhard Keller, {\em Deriving DG categories},
Annales de l'\'Ecole Normale Sup\'erieure,  {\bf 27} (1994) 63-102.

\bibitem{Kelleralgebraic}
Bernhard Keller, {\em On differential graded categories},
Proceedings of the International Congress of Mathematicians, Madrid, Spain,
2006, vol II, European Mathematical Society, 2006, pp. 151-190.

\bibitem{Krauseschemes}
Henning Krause, {\em The stable category of a Noetherian scheme}, Compositio Mathematica
{\bf 141} (2005) 1128-1162.

\bibitem{Matsumura} Hydeyuki Matsumura, {\sc Commutative Ring
Theory}. Cambridge University Press, Cambridge 1980.

\bibitem{May}
Peter May, {\em The axioms for triangulated categories}, manuscript\\
{\tt http://www.math.uchicago.edu/$\sim$may/MISCMaster.html}

\bibitem{Riedtmann}
Christine Riedtmann,
{\em Degenerations for representations of quivers with relations},
Annales de l'\'Ecole Normale Sup\'erieure {\bf 19} (1986) 275-301.

\bibitem{Verdier} Jean-Louis Verdier, {\sc Des cat\'egories d\'eriv\'ees des cat\'egories
ab\'eliennes}, Ast\'erisque {\bf 239} (1996).

\bibitem{YoshinoCM}
Yuji Yoshino, {\em On degeneration of Cohen-Macaulay modules},
Journal of Algebra {\bf 248} (2002) 272-290.

\bibitem{YoshinoModules}
Yuji Yoshino, {\em On degenerations of Modules}, Journal of Algebra {\bf 278} (2004) 217-226.

\bibitem{Yoshinostable}
Yuji Yoshino, {\em Stable degeneration of Cohen-Macaulay modules},
Journal of Algebra {\bf 332} (2011) 500-521.

\bibitem{reptheobuch}
Alexander Zimmermann, {\sc Representation Theory: A homological algebra point of view},
Springer Verlag, London 2014.

\bibitem{Osakaproceedings}
Alexander Zimmermann, {\em Remarks on a Categorical Definition of Degeneration in Triangulated Categories}, Proceedings of the 47th Symposion on Ring Theory and Representation Theory September 13-15, 2014. Saitama University, February 2015, pp 138-145.

\bibitem{Zwara}
Grzegorz Zwara, {\em Degenerations of finite dimensional modules are given by extensions},
Compositio Mathematica {\bf 121} (2000) 205-218.

\end{thebibliography}
\end{document}